\documentclass[12pt,a4paper]{article}
\usepackage{amssymb}

\usepackage{amsthm}
\newtheorem{definition}{Definition}[section]
\newtheorem{proposition}{Proposition}[section]
\newtheorem{theorem}{Theorem}[section]
\newtheorem{remark}{Remark}[section]
\newtheorem{corollary}{Corollary}[section]
\newtheorem{lemma}{Lemma}[section]

\begin{document}

\author{Nicoleta Aldea and Gheorghe Munteanu}
\title{Projectively related complex Finsler metrics}
\date{}
\maketitle

\begin{abstract}
In this paper we introduce in study the projectively related complex Finsler
metrics. We prove the complex versions of the Rapcs\'{a}k's theorem and
characterize the weakly K\"{a}hler and generalized Berwald projectively
related complex Finsler metrics. The complex version of Hilbert's Fourth
Problem is also pointed out. As an application, the projectiveness of a
complex Randers metric is described.
\end{abstract}

\begin{flushleft}
\strut \textbf{2010 Mathematics Subject Classification:} 53B40, 53C60.

\end{flushleft}

\section{Introduction}

\setcounter{equation}{0}The problem of projectively related real Finsler
metrics is quite old in geometry and its origin is formulated in Hilbert's
Fourth Problem: determine the metrics on an open subset in $\mathbf{R}^{n}$,
whose geodesics are straight lines. Two Finsler metrics, on a common
underlying manifold, are called projectively related if any geodesic of the
first is also geodesic for the second and the other way around.

The study of projective real Finsler spaces was initiated by L. Berwald,
\cite{Bw1,Bw2} , and his studies mainly concern the two dimensional Finsler
spaces. Further substantial contributions on this topic are from Rapcs\'{a}k
\cite{Ra}, Misra \cite{Mi} and, especially, from Z. Szabo \cite{Sz1} and M.
Matsumoto \cite{Map}. The problem of projective Finsler spaces is strongly
connected to projectively related sprays, as Z. Shen pointed out in \cite
{Sh1}. The topic of projective real Finsler spaces continues to be of
interest for special classes of metrics (\cite{B-M1,B-M2,C-ZS,Li,B-C}, etc.).

In complex geometry, T. Aikou studied in \cite{Ai} the projective flatness
of complex Finsler metrics by the projective flatness of Finsler connections.

Part of the general themes from projective real Finsler geometry can be
broached in complex Finsler geometry. However, there are meaningful
differences comparing to real reasonings, mainly on account of the fact that
the Chern-Finsler complex nonlinear connection (the main tool in this
geometry), generally does not derive from a spray. Another problem is that
in complex Finsler geometry, the notion of complex geodesic curve comports
two different nuances, one is in Abate-Patrizio's sense, (\cite{A-P}), and
the second is due to Royden, (\cite{Ro}). But, these notions don't differ
too much. Since a complex geodesic curve in Royden's sense assures that the
weakly K\"{a}hler condition is satisfied along the curve, we can state that
any complex geodesic curve in \cite{Ro}'s sense is a complex geodesic curve
in \cite{A-P}' s sense.

Our aim in the present paper is to study the projectively related complex
Finsler metrics $F$ and $\tilde{F}$ on the complex manifold $M$, using some
ideas from the real case. We have the canonical complex nonlinear connection
available, proven to be derived from a complex spray and hence it will
become an important tool in our approach. Also, in order to obtain a general
characterization of the projectively related complex Finsler metrics we
consider the complex geodesics in \cite{A-P}' s sense.

Subsequently, we have made an overview of the paper's content.

In \S 2, we recall some preliminary properties of the $n$ - dimensional
complex Finsler spaces.

In \S 3, we introduce the projectively related complex Finsler metrics and
then we find some necessary and sufficient condition of projectiveness,
(Theorem 3.1 and Corollary 3.1). For two projectively related complex
Finsler metrics we show that if one of these is weakly K\"{a}hler then, the
other must also be weakly K\"{a}hler, (Theorem 3.2). We prove some complex
versions of the Rapcs\'{a}k's theorem (Theorems 3.3, 3.4 and 3.5). Next, by
means of these theorems we are able to characterize the weakly K\"{a}hler
and generalized Berwald projectively related complex Finsler metrics,
(Theorem 3.6 and Corollary 3.2). Moreover, the complex version of Hilbert's
Fourth Problem is emphasized, (Theorem 3.7).

The last part of the paper (\S 4) is devoted to the projectiveness of the
complex Randers metric $\tilde{F}=\alpha +|\beta |.$ The necessary and
sufficient conditions in which the metrics $\tilde{F}$ and $\alpha $ are
projectively related are contained in Theorem 4.3. We prove that the complex
Randers metric $\tilde{F}=\alpha +|\beta |$ on a domain $D$ from $\mathbf{C}%
^n$ is projectively related to the complex Euclidean metric $F$ on $D$ if
and only if $\alpha $ is projectively related to the Euclidean metric $F$
and, $\tilde{F}$ is a complex Berwald metric, (Theorem 4.4).

\section{Preliminaries}

\setcounter{equation}{0}For the beginning we will make a survey of complex
Finsler geometry and we will set the basic notions and terminology. For more
see \cite{A-P,Mub, Al-Mu2}.

Let $M$ be a $n$ - dimensional complex manifold, $z=(z^k)_{k=\overline{1,n}}$
are complex coordinates in a local chart. The complexified of the real
tangent bundle $T_CM$ splits into the sum of holomorphic tangent bundle $%
T^{\prime }M$ and its conjugate $T^{\prime \prime }M$. The bundle $T^{\prime
}M$ is itself a complex manifold, and the local coordinates in a local chart
will be denoted by $u=(z^k,\eta ^k)_{k=\overline{1,n}}.$ These are changed
into $(z^{\prime k},\eta ^{\prime k})_{k=\overline{1,n}}$ by the rules $%
z^{\prime k}=z^{\prime k}(z)$ and $\eta ^{\prime k}=\frac{\partial z^{\prime
k}}{\partial z^l}\eta ^l.$

A \textit{complex Finsler space} is a pair $(M,F)$, where $F:T^{\prime
}M\rightarrow \mathbb{R}^{+}$ is a continuous function satisfying the
conditions:

\textit{i)} $L:=F^2$ is smooth on $\widetilde{T^{\prime }M}:=T^{\prime
}M\backslash \{0\};$

\textit{ii)} $F(z,\eta )\geq 0$, the equality holds if and only if $\eta =0;$

\textit{iii)} $F(z,\lambda \eta )=|\lambda |F(z,\eta )$ for $\forall \lambda
\in \mathbb{C}$;

\textit{iv)} the Hermitian matrix $\left( g_{i\bar{j}}(z,\eta )\right) $ is
positive definite, where $g_{i\bar{j}}:=\frac{\partial ^{2}L}{\partial \eta
^{i}\partial \bar{\eta}^{j}}$ is the fundamental metric tensor.
Equivalently, it means that the indicatrix is strongly pseudo-convex.

We say that a function $f$ on $T^{\prime }M$ is $(p,q)$ - \textit{homogeneous%
} with respect to $\eta $ iff $f(z,\lambda \eta )=\lambda ^p\bar{\lambda}%
^qf(z,\eta ),$ for any $\lambda \in \mathbb{C}$. For instance, $L:=F^2$ is a
$(1,1)$ - homogeneous function.

Roughly speaking, the geometry of a complex Finsler space consists of the
study of geometric objects of the complex manifold $T^{\prime }M$ endowed
with the Hermitian metric structure defined by $g_{i\bar{j}}.$ Thus, the
first step is to study sections of the complexified tangent bundle of $%
T^{\prime }M,$ which is decomposed into the sum $T_{C}(T^{\prime
}M)=T^{\prime }(T^{\prime }M)\oplus T^{\prime \prime }(T^{\prime }M)$.

Let $VT^{\prime }M\subset T^{\prime }(T^{\prime }M)$ be the vertical bundle,
locally spanned by $\{\frac{\partial }{\partial \eta ^{k}}\}$, and $%
VT^{\prime \prime }M$ be its conjugate. At this point, the idea of complex
nonlinear connection, briefly $(c.n.c.),$ is an instrument in
'linearization' of this geometry. A $(c.n.c.)$ is a supplementary complex
subbundle to $VT^{\prime }M$ in $T^{\prime }(T^{\prime }M)$, i.e. $T^{\prime
}(T^{\prime }M)=HT^{\prime }M\oplus VT^{\prime }M.$ The horizontal
distribution $H_{u}T^{\prime }M$ is locally spanned by $\{\frac{\delta }{%
\delta z^{k}}=\frac{\partial }{\partial z^{k}}-N_{k}^{j}\frac{\partial }{%
\partial \eta ^{j}}\}$, where $N_{k}^{j}(z,\eta )$ are the coefficients of
the $(c.n.c.)$. The pair $\{\delta _{k}:=\frac{\delta }{\delta z^{k}},\dot{%
\partial}_{k}:=\frac{\partial }{\partial \eta ^{k}}\}$ will be called the
adapted frame of the $(c.n.c.)$ which obey the change rules $\delta _{k}=%
\frac{\partial z^{\prime j}}{\partial z^{k}}\delta _{j}^{\prime }$ and $\dot{%
\partial}_{k}=\frac{\partial z^{\prime j}}{\partial z^{k}}\dot{\partial}%
_{j}^{\prime }.$ By conjugation everywhere we have obtained an adapted frame
$\{\delta _{\bar{k}},\dot{\partial}_{\bar{k}}\}$ on $T_{u}^{\prime \prime
}(T^{\prime }M).$ The dual adapted bases are $\{dz^{k},\delta \eta ^{k}\}$
and $\{d\bar{z}^{k},\delta \bar{\eta}^{k}\}.$

Certainly, a main problem in this geometry is to determine a $(c.n.c.)$
related only to the fundamental function of the complex Finsler space $%
(M,F). $

The next step is the action of a derivative law $D$ on the sections of $%
T_{C}(T^{\prime }M).$ A Hermitian connection $D$, of $(1,0)$ - type, which
satisfies in addition $D_{JX}Y=JD_{X}Y,$ for all $X$ horizontal vectors and $%
J$ the natural complex structure of the manifold, is called Chern-Finsler
connection (cf. \cite{A-P})$.$ Locally, it is given by the following
coefficients (cf. \cite{Mub}):
\begin{equation}
N_{j}^{i}:=g^{\overline{m}i}\frac{\partial g_{l\overline{m}}}{\partial z^{j}}%
\eta ^{l}=L_{lj}^{i}\eta ^{l}\;;\;L_{jk}^{i}:=g^{\overline{l}i}\delta
_{k}g_{j\overline{l}}\;\;;\;C_{jk}^{i}:=g^{\overline{l}i}\dot{\partial}%
_{k}g_{j\overline{l}},  \label{0}
\end{equation}
where here and further on $\delta _{k}$ is related to the Chern-Finsler $%
(c.n.c.)$ and $D_{\delta _{k}}\delta _{j}=L_{jk}^{i}\delta _{i},$ $D_{\dot{%
\partial}_{k}}\dot{\partial}_{j}=C_{jk}^{i}\dot{\partial}_{i}$.

Let us recall that in \cite{A-P}'s terminology, the complex Finsler space $%
(M,F)$ is \textit{strongly K\"{a}hler} iff $T_{jk}^i=0,$ \textit{K\"{a}hler}$%
\;$iff $T_{jk}^i\eta ^j=0$ and \textit{weakly K\"{a}hler }iff\textit{\ } $%
g_{i\overline{l}}T_{jk}^i\eta ^j\overline{\eta }^l=0,$ where $%
T_{jk}^i:=L_{jk}^i-L_{kj}^i.$ In \cite{C-S} it is proved that strongly
K\"{a}hler and K\"{a}hler notions actually coincide. We notice that in the
particular case of complex Finsler metrics which come from Hermitian metrics
on $M,$ so-called \textit{purely Hermitian metrics} in \cite{Mub}, i.e. $g_{i%
\overline{j}}=g_{i\overline{j}}(z)$, all these kinds of K\"{a}hler coincide.

The Chern-Finsler $(c.n.c.)$ generally, does not derive from a spray, but it
always determine a complex spray with the local coefficients $G^i=\frac
12N_j^i\eta ^j.$ Instead, $G^i$ induce a $(c.n.c.)$ by $\stackrel{c}{N_j^i}:=%
\dot{\partial}_jG^i$ called \textit{canonical} in \cite{Mub}, where it is
proved that it coincides with Chern-Finsler $(c.n.c.)$ if and only if the
complex Finsler metric is K\"{a}hler. Note that $2G^i=N_j^i\eta ^j=\stackrel{%
c}{N_j^i}\eta ^j,$ and so $\eta ^k\stackrel{c}{\delta _k}=\eta ^k\delta _k,$
where $\stackrel{c}{\delta _k}$ is related to canonical $(c.n.c.)$, i.e. $%
\stackrel{c}{\delta _k}:=\frac \partial {\partial z^k}-\stackrel{c}{N_k^j}%
\dot{\partial}_j.$ Additionally, in the K\"{a}hler case, we have $\stackrel{c%
}{\delta _k}=\delta _k$.

In \cite{Al-Mu2} we have proven that the complex Finsler space $(M,F)$ is
\textit{generalized Berwald }iff $\dot{\partial}_{\bar{h}}G^i=0$ and $(M,F)$
is a complex \textit{Berwald} space iff it is K\"{a}hler and generalized
Berwald.

\section{Projectively related complex Finsler metrics}

\setcounter{equation}{0}In Abate-Patrizio's sense, (\cite{A-P} p. 101), a
complex geodesic curve is given by $D_{T^h+\overline{T^h}}T^h=\theta
^{*}(T^h,\overline{T^h}),$ where $\theta ^{*}=g^{\bar{m}k}g_{i\bar{p}}(L_{%
\bar{j}\bar{m}}^{\bar{p}}-L_{\bar{m}\bar{j}}^{\bar{p}})dz^i\wedge d\bar{z}%
^j\otimes \delta _k,$ for which it is proven in \cite{Mub} that $\theta
^{*k}=2g^{\bar{j}k}\stackrel{c}{\delta _{\bar{j}}}L$ and $\theta ^{*i}$ is
vanishing if and only if the space is weakly K\"{a}hler. Thus, the equations
of a complex geodesic $z=z(s)$ of $(M,L)$, with $s$ a real parameter, in
\cite{A-P}' s sense can be rewritten as
\begin{equation}
\frac{d^2z^i}{ds^2}+2G^k(z(s),\frac{dz}{ds})=\theta ^{*i}(z(s),\frac{dz}{ds}%
)\;;\;i=\overline{1,n},  \label{1}
\end{equation}
where by $z^i(s),$ $i=\overline{1,n},$ we denote the coordinates along of
curve $z=z(s).$

We note that the functions $\theta ^{*i}$ are $(1,1)$ - homogeneous with
respect to $\eta ,$ i.e. $(\dot{\partial}_k\theta ^{*i})\eta ^k=\theta ^{*i}$
and $(\dot{\partial}_{\bar{k}}\theta ^{*i})\bar{\eta}^k=\theta ^{*i}.$

Let $\tilde{L}$ be another complex Finsler metric on the underlying manifold
$M.$

\begin{definition}
The complex Finsler metrics $L$ and $\tilde{L}$ on the manifold $M$, are
called projectively related if they have the same complex geodesics as point
sets.
\end{definition}

This means that for any complex geodesic $z=z(s)$ of $(M,L)$ there is a
transformation of its parameter $s$, $\tilde{s}=\tilde{s}(s),$ with $\frac{d%
\tilde{s}}{ds}>0,$ such that $z=z(\tilde{s}(s))$ is a geodesic of $(M,\tilde{%
L})$ and, conversely.

We suppose that $z=z(s)$ is a complex geodesic of $(M,L).$ Thus, it
satisfies (\ref{1}). Taking an arbitrary transformation of the parameter $%
t=t(s),$ with $\frac{dt}{ds}>0,$ the equations (\ref{1}) cannot in general
be preserved. Indeed, for the new parameter $t$ we have

\begin{center}
$\frac{dz^i}{ds}=\frac{dz^i}{dt}\frac{dt}{ds}\;;\;\frac{d^2z^i}{ds^2}=\frac{%
d^2z^i}{dt^2}\left( \frac{dt}{ds}\right) ^2+\frac{dz^i}{dt}\frac{d^2t}{ds^2}%
\;;\;\theta ^{*k}\left( z,\frac{dz}{ds}\right) =\left( \frac{dt}{ds}\right)
^2\theta ^{*k}(z,\frac{dz}{dt}).$
\end{center}

Then,

$\left[ \frac{d^2z^i}{dt^2}+2G^i(z,\frac{dz}{dt})-\theta ^{*i}(z,\frac{dz}{dt%
})\right] \left( \frac{dt}{ds}\right) ^2$

$=\frac{d^2z^i}{ds^2}-\frac{dz^i}{dt}\frac{d^2t}{ds^2}+2G^i(z,\frac{dz}{ds}%
)-\theta ^{*i}(z,\frac{dz}{ds})=-\frac{dz^i}{dt}\frac{d^2t}{ds^2}.$

Therefore, the equations (\ref{1}) in parameter $t$ are
\begin{equation}
\frac{d^{2}z^{i}}{dt^{2}}+2G^{i}(z(t),\frac{dz}{dt})=\theta ^{*i}(z(t),\frac{%
dz}{dt})-\frac{dz^{i}}{dt}\frac{d^{2}t}{ds^{2}}\frac{1}{\left( \frac{dt}{ds}%
\right) ^{2}}\;;\;i=\overline{1,n},  \label{2}
\end{equation}
which is equivalent to
\begin{equation}
\frac{\frac{d^{2}x^{i}}{dt^{2}}+2G^{i}(z,\frac{dz}{dt})-\theta ^{*i}(z,\frac{%
dz}{dt})}{\frac{dz^{i}}{dt}}=-\frac{d^{2}t}{ds^{2}}\frac{1}{\left( \frac{dt}{%
ds}\right) ^{2}}\;;\;i=\overline{1,n}.  \label{2'}
\end{equation}
We can rewrite (\ref{2'}), taking for $i$ two different values, as
\begin{equation}
\frac{\frac{d^{2}x^{i}}{dt^{2}}+2G^{j}(z,\frac{dz}{dt})-\theta ^{*j}(z,\frac{%
dz}{dt})}{\frac{dz^{j}}{dt}}=\frac{\frac{d^{2}z^{k}}{dt^{2}}+2G^{k}(z,\frac{%
dz}{dt})-\theta ^{*k}(z,\frac{dz}{dt})}{\frac{dz^{k}}{dt}}=-\frac{d^{2}t}{%
ds^{2}}\frac{1}{\left( \frac{dt}{ds}\right) ^{2}},  \label{3}
\end{equation}
for any $j,k=\overline{1,n}.$

Corresponding to the complex Finsler metric $\tilde{L}$ on the same manifold
$M,$ we have the spray coefficients $\tilde{G}^i$ and the functions $\tilde{%
\theta}^{*i}.$ If $L$ and $\tilde{L}$ are projectively related, then $z=z(%
\tilde{s})$ is a complex geodesic of $(M,\tilde{L}),$ where $\tilde{s}$ is
the parameter with respect to $\tilde{L}$ . Now, we assume that the same
parameter $t$ is transformed by $t=t(\tilde{s})$ and as above we obtain
\begin{equation}
\frac{\frac{d^2x^i}{dt^2}+2\tilde{G}^i(z,\frac{dz}{dt})-\tilde{\theta}%
^{*i}(z,\frac{dz}{dt})}{\frac{dz^i}{dt}}=-\frac{d^2t}{d\tilde{s}^2}\frac
1{\left( \frac{dt}{d\tilde{s}}\right) ^2}\;;\;i=\overline{1,n}.  \label{3'}
\end{equation}

The difference between (\ref{2'}) and (\ref{3'}) gives
\begin{equation}
2\tilde{G}^i(z,\frac{dz}{dt})-\tilde{\theta}^{*i}(z,\frac{dz}{dt})=2G^i(z,%
\frac{dz}{dt})-\theta ^{*i}(z,\frac{dz}{dt})+\left[ \frac{d^2t}{ds^2}\frac
1{\left( \frac{dt}{ds}\right) ^2}-\frac{d^2t}{d\tilde{s}^2}\frac 1{\left(
\frac{dt}{d\tilde{s}}\right) ^2}\right] \frac{dz^i}{dt}.  \label{3''}
\end{equation}

On the geodesic curves, it can be rewritten more generally as
\begin{equation}
2\tilde{G}^{i}(z,\frac{dz}{dt})-\tilde{\theta}^{*i}(z,\frac{dz}{dt}%
)=2G^{i}(z,\frac{dz}{dt})-\theta ^{*i}(z,\frac{dz}{dt})+2P(z,\frac{dz}{dt})%
\frac{dz^{i}}{dt},  \label{4}
\end{equation}
for any $i=\overline{1,n}$, where $P$ is a smooth function on $T^{\prime }M$
with complex values.

Denoting by $B^i:=\frac 12(\tilde{\theta}^{*i}-\theta ^{*i}),$ the
homogeneity properties of the functions $\tilde{\theta}^{*i}$ and $\theta
^{*i}$ give $(\dot{\partial}_kB^i)\eta ^k=B^i$ and $(\dot{\partial}_{\bar{k}%
}B^i)\bar{\eta}^k=B^i.$ Moreover the relations (\ref{4}) become
\begin{equation}
\tilde{G}^i=G^i+B^i+P\eta ^i.  \label{5}
\end{equation}

Now, we use their homogeneity properties, going from $\eta $ to $\lambda
\eta .$ Thus, differentiating in (\ref{5}) with respect to $\eta $ and $\bar{%
\eta}$ and then setting $\lambda =1$, we obtain
\begin{equation}
B^i=[(\dot{\partial}_kP)\eta ^k-P]\eta ^i\;\;\;\mbox{and}\;\;B^i=-(\dot{%
\partial}_{\bar{k}}P)\bar{\eta}^k\eta ^i  \label{6}
\end{equation}
and so,
\begin{equation}
(\dot{\partial}_kP)\eta ^k+(\dot{\partial}_{\bar{k}}P)\bar{\eta}^k=P,
\label{7}
\end{equation}
for any $i=\overline{1,n}$.

\begin{lemma}
Between the spray coefficients $\tilde{G}^i$ and $G^i$ of the metrics $L$
and $\tilde{L}$ on the manifold $M$ there are the relations $\tilde{G}%
^i=G^i+B^i+P\eta ^i,$ for any $i=\overline{1,n},$ where $P$ is a smooth
function on $T^{\prime }M$ with complex values, if and only if $\tilde{G}%
^i=G^i+(\dot{\partial}_kP)\eta ^k\eta ^i$, $B^i(z,\eta )=-(\dot{\partial}_{%
\bar{k}}P)\bar{\eta}^k\eta ^i$, for any $i=\overline{1,n},$ and $(\dot{%
\partial}_kP)\eta ^k+(\dot{\partial}_{\bar{k}}P)\bar{\eta}^k=P$.
\end{lemma}

From above considerations we obtain

\begin{lemma}
If the complex Finsler metrics $L$ and $\tilde{L}$ on the manifold $M$ are
projectively related, then there is a smooth function $P$ on $T^{\prime }M$
with complex values$,$ satisfying $(\dot{\partial}_kP)\eta ^k+(\dot{\partial}%
_{\bar{k}}P)\bar{\eta}^k=P,$ such that
\begin{equation}
\tilde{G}^i(z,\eta )=G^i(z,\eta )+(\dot{\partial}_kP)\eta ^k\eta ^i\;%
\mbox{and}\;B^i(z,\eta )=-(\dot{\partial}_{\bar{k}}P)\bar{\eta}^k\eta
^i\;;\;i=\overline{1,n}.  \label{8}
\end{equation}
\end{lemma}

\begin{remark}
We denote $S:=(\dot{\partial}_kP)\eta ^k$ and $Q:=-(\dot{\partial}_{\bar{k}%
}P)\bar{\eta}^k.$ The $(2,0)$-homogeneity with respect to $\eta $ of the
functions $\tilde{G}^i$ and $G^i$ implies the $(1,0)$-homogeneity of $S,$
and the $(1,1)$ - homogeneity of $B^i$ gives that $Q$ is $(0,1)$-homogeneous.
\end{remark}

Conversely, under assumption that $z=z(s)$ is a complex geodesic of $(M,L),$
we show that the complex Finsler metric $\tilde{L}$ with the spray
coefficients $\tilde{G}^i$ given by
\[
\tilde{G}^i=G^i+B^i+P\eta ^i,
\]
where $P$ is a smooth function on $T^{\prime }M$ with complex values, is
projectively related to $L,$ i.e. there is a parametrization $\tilde{s}=%
\tilde{s}(s),$ with $\frac{d\tilde{s}}{ds}>0,$ such that $z=z(\tilde{s}(s))$
is a geodesic of $(M,\tilde{L}).$

If there is a parametrization $\tilde{s}=\tilde{s}(s)$ then we have

$\frac{d^2z^i}{d\tilde{s}^2}=-2G^i(z,\frac{dz}{d\tilde{s}})+\theta ^{*i}(z,%
\frac{dz}{d\tilde{s}})-\frac{d^2\tilde{s}}{ds^2}\frac 1{\left( \frac{d\tilde{%
s}}{ds}\right) ^2}\frac{dz^i}{d\tilde{s}},$ for any for any $i=\overline{1,n}
$.

Now, using (\ref{8}), it results
\[
\frac{d^2z^i}{d\tilde{s}^2}=-2\tilde{G}^i(z,\frac{dz}{dt})+\tilde{\theta}%
^{*i}(z,\frac{dz}{d\tilde{s}})+\left( 2P(z,\frac{dz}{d\tilde{s}})-\frac{d^2%
\tilde{s}}{ds^2}\frac 1{\left( \frac{d\tilde{s}}{ds}\right) ^2}\right) \frac{%
dz^i}{d\tilde{s}}\;;\;i=\overline{1,n}.
\]
So, $z=z(\tilde{s}(s))$ is a geodesic of $(M,\tilde{L})$ if and only if
\begin{equation}
\left( 2P(z,\frac{dz}{d\tilde{s}})-\frac{d^2\tilde{s}}{ds^2}\frac 1{\left(
\frac{d\tilde{s}}{ds}\right) ^2}\right) \frac{dz^i}{d\tilde{s}}=0\;;\;i=%
\overline{1,n}.  \label{9}
\end{equation}
Supposing the complex geodesic curve is not a line, it results
\begin{equation}
2P(z,\frac{dz}{ds})\frac{d\tilde{s}}{ds}=\frac{d^2\tilde{s}}{ds^2}.
\label{10}
\end{equation}
Denoting by $u(s):=\frac{d\tilde{s}}{ds}$, we have $\frac{d^2\tilde{s}}{ds^2}%
=\frac{du}{ds}$ and so, $2P(z,\frac{dz}{ds})u=\frac{du}{ds}.$ We obtain $%
u=ae^{\int 2P(z,\frac{dz}{ds})ds}.$ From here, it results that there is
\[
\tilde{s}(s)=a\int e^{\int 2P(z,\frac{dz}{ds})ds}ds+b,
\]
where $a,b$ are arbitrary constants.

Corroborating all above results we have proven.

\begin{theorem}
Let $L$ and $\tilde{L}$ be complex Finsler metrics on the manifold $M$. Then
$L$ and $\tilde{L}$ are projectively related if and only if there is a
smooth function $P$ on $T^{\prime }M$ with complex values, such that
\begin{equation}
\tilde{G}^{i}=G^{i}+B^{i}+P\eta ^{i};\;i=\overline{1,n}.  \label{12}
\end{equation}
\end{theorem}

As a consequence of Lemma 3.1 we have the following.

\begin{corollary}
Let $L$ and $\tilde{L}$ be the complex Finsler metrics on the manifold $M$. $%
L$ and $\tilde{L}$ are projectively related if and only if there is a smooth
function $P$ on $T^{\prime }M$ with complex values, such that $\tilde{G}%
^{i}=G^{i}+(\dot{\partial}_{k}P)\eta ^{k}\eta ^{i}$ , $B^{i}(z,\eta )=-(\dot{%
\partial}_{\bar{k}}P)\bar{\eta}^{k}\eta ^{i}$, for any $i=\overline{1,n},$
and $(\dot{\partial}_{k}P)\eta ^{k}+(\dot{\partial}_{\bar{k}}P)\bar{\eta}%
^{k}=P$.
\end{corollary}

The relations (\ref{12}) between the spray coefficients $\tilde{G}^i$ and $%
G^i$ of the projectively related complex Finsler metrics $L$ and $\tilde{L}$
will be called \textit{projective change.}

\begin{theorem}
Let $L$ and $\tilde{L}$ be two complex Finsler metrics on the manifold $M,$
which are projectively related. Then, $L$ is weakly K\"{a}hler if and only
if $\tilde{L}$ is also weakly K\"{a}hler. In this case, the projective
change is $\tilde{G}^{i}=G^{i}+P\eta ^{i}$, where $P$ is a $(1,0)$ -
homogeneous function.
\end{theorem}

\begin{proof} We assume that $\tilde{G}^i=G^i+(\dot{\partial}_kP)\eta
^k\eta ^i$, $B^i=\frac 12(\tilde{\theta}^{*i}-\theta ^{*i})=-(\dot{\partial}%
_{\bar{k}}P)\bar{\eta}^k\eta ^i$ and $(\dot{\partial}_kP)\eta ^k+(\dot{%
\partial}_{\bar{k}}P)\bar{\eta}^k=P.$

If $L$ is weakly K\"{a}hler then $\theta ^{*i}=0$ and so, $\tilde{\theta}%
^{*i}=-2(\dot{\partial}_{\bar{k}}P)\bar{\eta}^k\eta ^i$ which contracted by $%
\tilde{g}_{i\bar{r}}\bar{\eta}^r=\dot{\partial}_l\tilde{L},$ gives $\tilde{%
\theta}^{*i}\tilde{g}_{i\bar{r}}\bar{\eta}^r=-2(\dot{\partial}_{\bar{k}}P)%
\bar{\eta}^k\tilde{L}.$

But, $\tilde{\theta}^{*i}\tilde{g}_{i\bar{r}}\bar{%
\eta}^r=0.$ Thus, $(\dot{\partial}_{\bar{k}}P)\bar{\eta}^k=0$, which implies $%
\tilde{\theta}^{*i}=0,$ i.e. $\tilde{L}$ is weakly K\"{a}hler and $P=(\dot{%
\partial}_kP)\eta ^k.$ So that, we obtain $\tilde{G}^i=G^i+P\eta ^i.$ The
converse implication results immediately by the same way.
\end{proof}

\begin{lemma}
Let $(M,L)$ be a complex Finsler space and $\tilde{L}$ a complex Finsler
metric on $M.$ The spray coefficients $\tilde{G}^{i}$ and $G^{i}$ of the
metrics $L$ and $\tilde{L}$ satisfy
\begin{equation}
\tilde{G}^{i}=G^{i}+\frac{1}{2}\tilde{g}^{\bar{r}i}\left( \dot{\partial}_{%
\bar{r}}(\delta _{k}\tilde{L})\eta ^{k}+2(\dot{\partial}_{\bar{r}}G^{l})(%
\dot{\partial}_{l}\tilde{L})\right) ;\;i=\overline{1,n}.  \label{11}
\end{equation}
\end{lemma}

\begin{proof} Having $\stackrel{c}{\delta _k}\tilde{L}=\frac{\partial
\tilde{L}}{\partial z^k}-\stackrel{c}{N_k^l}(\dot{\partial}_l\tilde{L}),$ by
a direct computation we obtain

$\dot{\partial}_{\bar{r}}(\stackrel{c}{\delta _k}\tilde{L})=\dot{\partial}_{%
\bar{r}}\left( \frac{\partial \tilde{L}}{\partial z^k}-\stackrel{c}{N_k^l}(%
\dot{\partial}_l\tilde{L})\right) =\frac{\partial ^2\tilde{L}}{\partial
z^k\partial \bar{\eta}^r}-(\dot{\partial}_{\bar{r}}\stackrel{c}{N_k^l})(\dot{%
\partial}_l\tilde{L})-\stackrel{c}{N_k^l}\tilde{g}_{l\bar{r}},$ which
contracted with $\tilde{g}^{\bar{r}i}\eta ^k$, and taking into account $\eta^k \stackrel{c}{\delta _k}=\eta^k \delta_k$, implies that

$\tilde{g}^{\bar{r}i}\dot{\partial}_{\bar{r}}(\stackrel{c}{\delta _k}\tilde{L%
})\eta ^k=\tilde{g}^{\bar{r}i}\dot{\partial}_{\bar{r}}(\delta _k\tilde{L}%
)\eta ^k=2\tilde{G}^i-2\tilde{g}^{\bar{r}i}(\dot{\partial}_{\bar{r}}G^l)(%
\dot{\partial}_l\tilde{L})-2G^i$ and so (\ref{11}) is justified.
\end{proof}

Next, we prove some complex versions of the Rapcs\'{a}k's theorem.

\begin{theorem}
Let $L$ and $\tilde{L}$ be complex Finsler metrics on the manifold $M$.
Then, $L$ and $\tilde{L}$ are projectively related if and only if
\begin{equation}
\frac{1}{2}\left( \dot{\partial}_{\bar{r}}(\delta _{k}\tilde{L})\eta ^{k}+2(%
\dot{\partial}_{\bar{r}}G^{l})(\dot{\partial}_{l}\tilde{L})\right) =P(\dot{%
\partial}_{\bar{r}}\tilde{L})+B^{i}\tilde{g}_{i\bar{r}}\;;\;r=\overline{1,n},
\label{13}
\end{equation}
with $\displaystyle P=\frac{1}{2\tilde{L}}[(\delta _{k}\tilde{L})\eta
^{k}+\theta ^{*i}(\dot{\partial}_{i}\tilde{L})].$
\end{theorem}

\begin{proof} We assume that $L$ and $\tilde{L}$ are projectively
related. Then, by Theorem 3.1 and (\ref{11}) we have
\begin{equation}
B^i+P\eta ^i=\frac 12\tilde{g}^{\bar{r}i}\left( \dot{\partial}_{\bar{r}%
}(\delta _k\tilde{L})\eta ^k+2(\dot{\partial}_{\bar{r}}G^l)(\dot{\partial}_l%
\tilde{L})\right) ;\;i=\overline{1,n}.  \label{14}
\end{equation}
First, if these relations are contracted by $\tilde{g}_{i\bar{m}}\bar{\eta}%
^m$, we get
\[
-\frac 12\theta ^{*i}(\dot{\partial}_i\tilde{L})+P\tilde{L}=\frac 12\dot{%
\partial}_{\bar{m}}(\delta _k\tilde{L})\eta ^k\bar{\eta}^m+(\dot{\partial}_{%
\bar{m}}G^l)\bar{\eta}^m(\dot{\partial}_l\tilde{L}),
\]
because $B^i\tilde{g}_{i\bar{m}}\bar{\eta}^m=-\frac 12\theta ^{*i}(\dot{%
\partial}_i\tilde{L}).$ But, the $(2,0)-$ homogeneity of the functions $G^i$ leads to $(\dot{\partial}_{\bar{m}}G^l)\bar{\eta}^m=0$
and $\dot{\partial}_{\bar{m}}(\delta _k\tilde{L})\eta ^k\bar{\eta}^m=(\delta
_k\tilde{L})\eta ^k.$ Thus, $P=\frac 1{2\tilde{L}}[(\delta _k\tilde{L})\eta
^k+\theta ^{*i}(\dot{\partial}_i\tilde{L})]$. Second, contracting into (%
\ref{14}) only by $\tilde{g}_{i\bar{m}},$ we obtain (\ref{13}).

Conversely, plugging the formulas (\ref{13}) into (\ref{11}), it results (%
\ref{12}) with $P=\frac 1{2\tilde{L}}[(\delta _k\tilde{L})\eta ^k+\theta
^{*i}(\dot{\partial}_i\tilde{L})]$, i.e. $L$ and $\tilde{L}$ are
projectively related.
\end{proof}

\begin{theorem}
Let $L$ and $\tilde{L}$ be the complex Finsler metrics on the manifold $M$.
Then, $L$ and $\tilde{L}$ are projectively related if and only if
\begin{eqnarray}
\dot{\partial}_{\bar{r}}(\delta _{k}\tilde{L})\eta ^{k}+2(\dot{\partial}_{%
\bar{r}}G^{l})(\dot{\partial}_{l}\tilde{L}) &=&\frac{1}{\tilde{L}}(\delta
_{k}\tilde{L})\eta ^{k}(\dot{\partial}_{\bar{r}}\tilde{L});  \label{14'} \\
B^{r} &=&-\frac{1}{2\tilde{L}}\theta ^{*l}(\dot{\partial}_{l}\tilde{L})\eta
^{r};\;r=\overline{1,n};  \nonumber \\
P &=&\frac{1}{2\tilde{L}}[(\delta _{k}\tilde{L})\eta ^{k}+\theta ^{*i}(\dot{%
\partial}_{i}\tilde{L})].  \nonumber
\end{eqnarray}
Moreover, the projective change is $\tilde{G}^{i}=G^{i}+\frac{1}{2\tilde{L}}%
(\delta _{k}\tilde{L})\eta ^{k}\eta ^{i}.$
\end{theorem}

\begin{proof} By Corollary 3.1, if $L$ and $\tilde{L}$ are
projectively related, then there is a smooth function $P$ on $T^{\prime }M$
with complex values, such that $\tilde{G}^i=G^i+(\dot{\partial}_kP)\eta
^k\eta ^i$, $B^i=-(\dot{\partial}_{\bar{k}}P)\bar{\eta}^k\eta ^i$, for any $%
i=\overline{1,n},$ and $(\dot{\partial}_kP)\eta ^k+(\dot{\partial}_{\bar{k}%
}P)\bar{\eta}^k=P$. Using (\ref{11}) it results
\begin{equation}
(\dot{\partial}_kP)\eta ^k\eta ^i=\frac 12\tilde{g}^{\bar{r}i}\left( \dot{%
\partial}_{\bar{r}}(\delta _k\tilde{L})\eta ^k+2(\dot{\partial}_{\bar{r}%
}G^l)(\dot{\partial}_l\tilde{L})\right) ;\;i=\overline{1,n},  \label{14''}
\end{equation}
which contracted firstly by $\tilde{g}_{i\bar{m}}$ and secondly by $\tilde{g}%
_{i\bar{m}}\bar{\eta}^m$ give

$\dot{\partial}_{\bar{r}}(\delta _k\tilde{L})\eta ^k+2(\dot{\partial}_{\bar{r%
}}G^l)(\dot{\partial}_l\tilde{L})=2(\dot{\partial}_kP)\eta ^k(\dot{\partial}%
_{\bar{r}}\tilde{L})$ and $(\dot{\partial}_kP)\eta ^k=\frac 1{2\tilde{L}%
}(\delta _k\tilde{L})\eta ^k$ respectively, since $\delta_k \tilde{L}$ is $(1,1)-$homogeneous.

Now, contracting $B^i=-(\dot{\partial}_{\bar{k}}P)\bar{\eta}^k\eta ^i$ with $%
\tilde{g}_{i\bar{m}}\bar{\eta}^m$, it leads to $(\dot{\partial}_{\bar{k}}P)%
\bar{\eta}^k=\frac 1{2\tilde{L}}\theta ^{*i}(\dot{\partial}_i\tilde{L}),$
because $B^i$ $\tilde{g}_{i\bar{m}}\bar{\eta}^m=-\frac 12\theta ^{*i}(\dot{%
\partial}_i\tilde{L}).$ Adding the last two relations obtained, it results $%
P=\frac 1{2\tilde{L}}[(\delta _k\tilde{L})\eta ^k+\theta ^{*i}(\dot{\partial}%
_i\tilde{L})].$

Conversely, replacing the first condition of (\ref{14'}) into (\ref{11}) we
obtain $\tilde{G}^i=G^i+S\eta ^i$, where $S:=\frac 1{2\tilde{L}}(\delta _k%
\tilde{L})\eta ^k.$

Now, having $P=\frac 1{2\tilde{L}}[(\delta _k\tilde{L})\eta ^k+\theta ^{*i}(%
\dot{\partial}_i\tilde{L})],$ we obtain

\begin{center}
$(\dot{\partial}_kP)\eta ^k=\frac 1{2\tilde{L}}(\delta _k\tilde{L})\eta ^k=S$
and $(\dot{\partial}_{\bar{k}}P)\bar{\eta}^k=\frac 1{2\tilde{L}}\theta ^{*i}(%
\dot{\partial}_i\tilde{L}).$
\end{center}

Thus, these lead to $\tilde{G}^i=G^i+(\dot{\partial}_kP)\eta ^k\eta ^i$, $%
B^i=-(\dot{\partial}_{\bar{k}}P)\bar{\eta}^k\eta ^i$ and $(\dot{\partial}%
_kP)\eta ^k+(\dot{\partial}_{\bar{k}}P)\bar{\eta}^k=P.$
\end{proof}

Plugging $\tilde{L}=\tilde{F}^2$ into (\ref{14'}) we have proven another
equivalent complex version of Rapcs\'{a}k's theorem.

\begin{theorem}
Let $F$ and $\tilde{F}$ be the complex Finsler metrics on the manifold $M$.
Then, $F$ and $\tilde{F}$ are projectively related if and only if
\begin{eqnarray}
\dot{\partial}_{\bar{r}}(\delta _{k}\tilde{F})\eta ^{k}+2(\dot{\partial}_{%
\bar{r}}G^{l})(\dot{\partial}_{l}\tilde{F}) &=&\frac{1}{\tilde{F}}(\delta
_{k}\tilde{F})\eta ^{k}(\dot{\partial}_{\bar{r}}\tilde{F})\;;  \label{15} \\
B^{r} &=&-\frac{1}{\tilde{F}}\theta ^{*l}(\dot{\partial}_{l}\tilde{F})\eta
^{r};\;r=\overline{1,n};  \nonumber \\
P &=&\frac{1}{\tilde{F}}[(\delta _{k}\tilde{F})\eta ^{k}+\theta ^{*i}(\dot{%
\partial}_{i}\tilde{F})].  \nonumber
\end{eqnarray}
Moreover, the projective change is $\tilde{G}^{i}=G^{i}+\frac{1}{\tilde{F}}%
(\delta _{k}\tilde{F})\eta ^{k}\eta ^{i}.$
\end{theorem}

\begin{theorem}
Let $L$ be a weakly K\"{a}hler complex Finsler metric on the manifold $M$
and $\tilde{L}$ another complex Finsler metric on $M.$ Then, $L$ and $\tilde{%
L}$ are projectively related if and only if $\tilde{L}$ is weakly K\"{a}hler
and
\begin{eqnarray}
\dot{\partial}_{\bar{r}}(\delta _{k}\tilde{L})\eta ^{k}+2(\dot{\partial}_{%
\bar{r}}G^{l})(\dot{\partial}_{l}\tilde{L}) &=&2P(\dot{\partial}_{\bar{r}}%
\tilde{L})\;;\;r=\overline{1,n},  \label{16} \\
P &=&\frac{1}{2\tilde{L}}(\delta _{k}\tilde{L})\eta ^{k}.  \nonumber
\end{eqnarray}
Moreover, the projective change is $\tilde{G}^{i}=G^{i}+P\eta ^{i}$ and $P$
is $(1,0)$ - homogeneous.
\end{theorem}

\begin{proof} Having in mind the Theorems 3.2 and 3.4 the direct
implication is obvious. For the converse, we have $B^i=\theta ^{*i}=\tilde{%
\theta}^{*i}=0,$ because $L$ and $\tilde{L}$ are weakly K\"{a}hler, which
together with (\ref{16}) are sufficient conditions for the projectiveness of
the metrics $L$ and $\tilde{L}$. Now, plugging (\ref{16}) into (\ref{11})
it results $\tilde{G}^i=G^i+P\eta ^i$ and the $(1,0)-$homogeneity of $P$.
\end{proof}

Let us pay more attention to Theorem 3.5. As its consequence, we have.

\begin{corollary}
Let $F$ be a generalized Berwald metric on the manifold $M$ and $\tilde{F}$
another complex Finsler metric on $M.$ Then, $F$ and $\tilde{F}$ are
projective if and only if
\begin{eqnarray}
\dot{\partial}_{\bar{r}}(\delta _{k}\tilde{F})\eta ^{k} &=&\frac{1}{\tilde{F}%
}(\delta _{k}\tilde{F})\eta ^{k}(\dot{\partial}_{\bar{r}}\tilde{F}%
)\;;\;B^{r}=-\frac{1}{\tilde{F}}\theta ^{*l}(\dot{\partial}_{l}\tilde{F}%
)\eta ^{r};  \label{19} \\
P &=&\frac{1}{\tilde{F}}[(\delta _{k}\tilde{F})\eta ^{k}+\theta ^{*i}(\dot{%
\partial}_{i}\tilde{F})],  \nonumber
\end{eqnarray}
for any $r=\overline{1,n}.$ Moreover, the projective change is $\tilde{G}%
^{i}=G^{i}+\frac{1}{\tilde{F}}(\delta _{k}\tilde{F})\eta ^{k}\eta ^{i}$ and $%
\tilde{F}$ is also generalized Berwald.
\end{corollary}

\begin{proof} The equivalence results by Theorem 3.5 in which $\dot{%
\partial}_{\bar{r}}G^l=0$, because $F$ is a generalized Berwald metric. In
order to show that $\tilde{F}$ is generalized Berwald, we compute

$\dot{\partial}_{\bar{r}}[\frac 1{\tilde{F}}(\delta _k\tilde{F})\eta
^k]=-\frac 1{\tilde{F}^2}(\dot{\partial}_{\bar{r}}\tilde{F})(\delta _k\tilde{%
F})\eta ^k+\frac 1{\tilde{F}}\dot{\partial}_{\bar{r}}(\delta _k\tilde{F}%
)\eta ^k=0,$ by using the first identity from (\ref{19}). Now,
differentiating the projective change $\tilde{G}^i=G^i+\frac 1{\tilde{F%
}}(\delta _k\tilde{F})\eta ^k\eta ^i$ with respect to $\bar{\eta}^r$ it
results $\dot{\partial}_{\bar{r}}\tilde{G}^l=0,$ i.e. $\tilde{F}$ is
generalized Berwald.
\end{proof}

In particular, if $F$ is a K\"{a}hler metric, then Theorem 3.6 and Corollary
3.2 imply

\begin{corollary}
Let $F$ be a complex Berwald metric on the manifold $M$ and $\tilde{F}$
another complex Finsler metric on $M.$ Then, $F$ and $\tilde{F}$ are
projectively related if and only if $\tilde{F}$ is weakly K\"{a}hler and
\begin{equation}
\dot{\partial}_{\bar{r}}(\delta _{k}\tilde{F})\eta ^{k}=P(\dot{\partial}_{%
\bar{r}}\tilde{F})\;;\;r=\overline{1,n}\;;\;P=\frac{1}{\tilde{F}}(\delta _{k}%
\tilde{F})\eta ^{k}.  \label{19''}
\end{equation}
Moreover, the projective change is $\tilde{G}^{i}=G^{i}+P\eta ^{i}$ and $%
\tilde{F}$ is generalized Berwald.
\end{corollary}

\begin{proposition}
Let $F$ and $\tilde{F}$ be two projectively related complex Finsler metrics
on the manifold $M$. If $P$ is $(1,0)$ - homogeneous with respect to $\eta $
and $F$ is generalized Berwald, then $P$ is holomorphic with respect to $%
\eta $.
\end{proposition}

\begin{proof} We have $\tilde{G}^i=G^i+B^i+P\eta ^i$, with $P$ homogenous
of $(1,0)$ - degree. This implies $B^i=0$ and so, by Corollary 3.2, $\theta
^{*l}(\dot{\partial}_l\tilde{F})=0.$ So that, $P=\frac 1{\tilde{F}}({\delta _k}\tilde{F})\eta ^k\eta ^i$ and, it has the property $\dot{\partial}_{%
\bar{r}}P=0.$
\end{proof}

\begin{proposition}
Let $F$ and $\tilde{F}$ be two projectively related complex Finsler metrics
on the manifold $M$. If $P$ is $(0,1)$ - homogeneous with respect to $\eta $
and $F$ is generalized Berwald, then $B^{i}=-P\eta ^{i}$, for any $i=%
\overline{1,n}$, and the projective change is $\tilde{G}^{i}=G^{i}$.
\end{proposition}

\begin{proof} Let be the projective change $\tilde{G}^i=G^i+B^i+P\eta ^i$, with $P$
homogenous of $(0,1)$ - degree. Since  $\tilde{G}^i$ and $G^i$ are $(2,0)$ - homogeneous and $B^i\;, P\eta^i$ are $(1,1)$ - homogeneous, it follows  $\tilde{G}^i=G^i$
and $B^i=-P\eta ^i$.
\end{proof}

Further on, the complex version of the \textit{Hilbert's Fourth Problem} is
approached.

\begin{theorem}
Let $L$ be complex Euclidean metric on a domain $D$ from $\mathbf{C}^{n}$
and $\tilde{L}$ another complex Finsler metric on $D.$ Then, $L$ and $\tilde{%
L}$ are projectively related if and only if $\tilde{L}$ is weakly K\"{a}hler
and
\begin{equation}
\tilde{G}^{i}=\frac{1}{2\tilde{L}}\frac{\partial \tilde{L}}{\partial z^{k}}%
\eta ^{k}\eta ^{i}\;;\;i=\overline{1,n}\;.  \label{17}
\end{equation}

Moreover, $\tilde{L}$ is generalized Berwald.
\end{theorem}

\begin{proof} The complex Euclidean metric $L:=|\eta |^2={\sum }%
_{k=1}^n\eta ^k\bar{\eta}^k$ is K\"{a}hler with the local spray coefficients $%
G^i=0 $, for any $i=\overline{1,n}$. By these assumptions, the conditions (%
\ref{16}) can be rewritten as
\begin{equation}
\dot{\partial}_{\bar{r}}(\frac{\partial \tilde{L}}{\partial z^k})\eta ^k=2P(%
\dot{\partial}_{\bar{r}}\tilde{L}),  \label{17'}
\end{equation}
for any $r=\overline{1,n}\;,$ where$\;P=\frac 1{2\tilde{L}}\frac{\partial
\tilde{L}}{\partial z^k}\eta ^k.$ Further on, by contraction in (\ref{17'})
with $\tilde{g}^{\bar{r}i}$ and since $\tilde G^i=\frac{1}{2}\tilde g^{\bar ri}\dot{\partial}_{\bar{r}}(\frac{\partial \tilde{L}}{\partial z^k})\eta ^k$,
using again (\ref{17'}) it follows that $\tilde G^i=P\eta^i$ which is (\ref{17}).
The converse is obvious.
\end{proof}

Taking $\tilde{L}=\tilde{F}^2$ into (\ref{17}), it becomes
\begin{equation}
\tilde{G}^i=\frac 1{\tilde{F}}\frac{\partial \tilde{F}}{\partial z^k}\eta
^k\eta ^i\;;\;i=\overline{1,n}\;.  \label{17''}
\end{equation}

Some examples of complex Finsler metrics which are projectively related to
the complex Euclidean metric are given by the following purely Hermitian
metrics defined over the disk $\Delta _{r}^{n}=\left\{ z\in \mathbf{C}%
^{n},\;|z|<r,\;\;r:=\sqrt{\frac{1}{|\varepsilon |}}\right\} $ $\subset
\mathbf{C}^{n}:$
\begin{equation}
\tilde{L}(z,\eta ):=\frac{|\eta |^{2}+\varepsilon \left( |z|^{2}|\eta
|^{2}-\left| <z,\eta >\right| ^{2}\right) }{(1+\varepsilon |z|^{2})^{2}}%
\;;\;\varepsilon <0,  \label{18}
\end{equation}
where $|z|^{2}:=$ ${\sum }_{k=1}^{n}z^{k}\bar{z}^{k},$ $<z,\eta >:={\sum }%
_{k=1}^{n}z^{k}\bar{\eta}^{k}$ and $\left| <z,\eta >\right| ^{2}=<z,\eta >%
\overline{<z,\eta >}$ $.$ They are K\"{a}hler and in particular, for $%
\varepsilon =-1$ we obtain the \textit{Bergman metric} on the unit disk $%
\Delta ^{n}:=\Delta _{1}^{n}$. Their geodesics are segments of straight
lines.

\section{Projectiveness of a complex Randers metric}

\setcounter{equation}{0}We consider $\beta (z,\eta ):=b_i(z)\eta ^i$ a
differential ($1,0)$ - form and $\alpha (z,\eta )$ :$=\sqrt{a_{i\bar{j}%
}(z)\eta ^i\bar{\eta}^j}$ a purely Hermitian metric on the manifold $M.$ By
these objects we have defined the complex Randers metric $\tilde{F}=\alpha
+|\beta |$ on $T^{\prime }M$ with
\begin{eqnarray*}
\frac{\partial \alpha }{\partial \eta ^i} &=&\frac 1{2\alpha }l_i\;\;;\;\;%
\frac{\partial |\beta |}{\partial \eta ^i}=\frac{\bar{\beta}}{2|\beta |}%
b_i\;;\;\;\tilde{\eta}_i:=\frac{\partial \tilde{L}}{\partial \eta ^i}=\frac{%
\tilde{F}}\alpha l_i+\frac{\tilde{F}\bar{\beta}}{|\beta |}b_i, \\
\tilde{G}^i &=&\stackrel{a}{G^i}+\frac 1{2\gamma }\left( l_{\bar{r}}\frac{%
\partial \bar{b}^r}{\partial z^j}-\frac{\beta ^2}{|\beta |^2}\frac{\partial
b_{\bar{r}}}{\partial z^j}\bar{\eta}^r\right) \xi ^i\eta ^j+\frac \beta
{4|\beta |}k^{\overline{r}i}\frac{\partial b_{\bar{r}}}{\partial z^j}\eta ^j
\\
l_i :&=&a_{i\bar{j}}\bar{\eta}^j\;;\;\;b^i:=a^{\bar{j}i}b_{\bar{j}%
}\;\;;\;\;||b||^2:=a^{\bar{j}i}b_ib_{\bar{j}}\;\;;\;\;b^{\bar{\imath}}:=\bar{%
b}^i,
\end{eqnarray*}
where $\stackrel{a}{G^i}=\frac 12\stackrel{a}{N_j^i}\eta ^j$ are the spray
coefficients of the purely Hermitian metric $\alpha $ and $\gamma :=\tilde{L}%
+\alpha ^2(||b||^2-1)$, $\xi ^i:=\bar{\beta}\eta ^i+\alpha ^2b^i$, $k^{\bar{r%
}i}:=2\alpha a^{\bar{j}i}+\frac{2(\alpha ||b||^2+2|\beta |)}\gamma \eta ^i%
\bar{\eta}^r-\frac{2\alpha ^3}\gamma b^i\bar{b}^r-\frac{2\alpha }\gamma (%
\bar{\beta}\eta ^i\bar{b}^r+\beta b^i\bar{\eta}^r)$.

Moreover, the complex Randers metric $\tilde{F}$ is weakly K\"{a}hler if and
only if
\[
\frac{\alpha ^2|\beta |}{\gamma \delta }\left[ \beta \frac{\alpha
||b||^2+|\beta |}{|\beta |}\frac{\partial b_{\bar{m}}}{\partial z^r}\bar{\eta%
}^m+\bar{\beta}\left( \frac{\partial b_r}{\partial z^l}-b^{\bar{m}}\frac{%
\partial a_{l\bar{m}}}{\partial z^r}\right) \eta ^l-\alpha |\beta |b^{\bar{m}%
}\frac{\partial b_{\bar{m}}}{\partial z^r}\right] \eta ^rC_k
\]
\begin{equation}
-\left( \alpha \bar{\beta}F_{kl}+\alpha b_l\frac{\partial b_{\bar{r}}}{%
\partial z^k}\bar{\eta}^r+2|\beta |a_{l\bar{r}}\Gamma _{\bar{j}k}^{\bar{r}}%
\bar{\eta}^j\right) \eta ^l+\alpha b_k\frac{\partial b_{\bar{m}}}{\partial
z^r}\bar{\eta}^m\eta ^r=0,  \label{19'}
\end{equation}
where $C_j:=\delta \left( \frac 1{\alpha ^2}l_j-\frac{\bar{\beta}}{|\beta |^2%
}b_j\right) ,$ with $\delta :=\frac{\alpha ^2||b||^2-|\beta |^2}{2\gamma }-%
\frac{n|\beta |}{2F},$ $F_{il}:=\frac{\partial b_l}{\partial z^i}-\frac{%
\partial bi}{\partial z^l}$ and $\Gamma _{\bar{j}i}^{\bar{r}}:=\frac 12a^{%
\bar{r}k}(\frac{\partial a_{k\bar{j}}}{\partial z^i}$ $-\frac{\partial a_{i%
\bar{j}}}{\partial z^k})$ . For more details see \cite{Al-Mu1}.

\begin{theorem}
(\cite{Al-Mu2}) Let $(M,\tilde{F})$ be a connected complex Randers space.
Then, it is a generalized Berwald space if and only if $(\bar{\beta}l_{\bar{r%
}}\frac{\partial \bar{b}^{r}}{\partial z^{j}}+\beta \frac{\partial b_{\bar{r}%
}}{\partial z^{j}}\bar{\eta}^{r})\eta ^{j}=0.$
\end{theorem}

\begin{theorem}
(\cite{Al-Mu2}) Let $(M,\tilde{F})$ be a connected complex Randers space.
Then, it$\ $is a complex Berwald space if and only if it is both generalized
Berwald and weakly K\"{a}hler. Moreover, $\alpha $ is K\"{a}hler and $\tilde{%
N}_{j}^{i}=\stackrel{a}{N_{j}^{i}}.$
\end{theorem}

First, our aim is to determine the necessary and sufficient conditions in
which the complex Randers metric $\tilde{F}$ is projectively related to the
Hermitian metric $\alpha .$ A simple computation shows that,

\begin{equation}
(\delta _k\tilde{F})\eta ^k=(\delta _k|\beta |)\eta ^k=\frac 1{2|\beta |}(%
\bar{\beta}l_{\bar{r}}\frac{\partial \bar{b}^r}{\partial z^k}+\beta \frac{%
\partial b_{\bar{r}}}{\partial z^k}\bar{\eta}^r)\eta ^k,  \label{20}
\end{equation}
because $(\delta _k\alpha )\eta ^k=0$ and
\begin{equation}
\theta ^{*i}(\dot{\partial}_i\tilde{F})=-\frac{\bar{\beta}}{2|\beta |}\Gamma
_{i\bar{j}}^kb_k\eta ^i\bar{\eta}^j.  \label{21}
\end{equation}

Taking into account Theorem 4.1 we have proven

\begin{lemma}
Let $(M,\tilde{F})$ be a connected complex Randers space. Then, $(M,\tilde{F}%
)$ is a generalized Berwald space if and only if $(\delta _{k}|\beta |)\eta
^{k}=0.$
\end{lemma}

\begin{theorem}
Let $(M,\tilde{F})$ be a connected complex Randers space. Then,

i) $\alpha $ and $\tilde{F}$ are projectively related if and only if $\tilde{%
F}$ is generalized Berwald and $B^{i}=-P\eta ^{i},$ for any $i=\overline{1,n}%
,$ where $P=-\frac{\bar{\beta}}{2\tilde{F}|\beta |}\Gamma _{i\bar{j}%
}^{k}b_{k}\eta ^{i}\bar{\eta}^{j};$

ii) $\alpha $ is K\"{a}hler and $\alpha $ is projectively related to $\tilde{%
F}$ if and only if $\tilde{F}$ is a complex Berwald metric.

In these cases, the projective change is $\tilde{G}^{i}=\stackrel{a}{G^{i}}.$
\end{theorem}

\begin{proof} We first prove i). The purely Hermitian property of the metric $%
\alpha $ implies that it is generalized Berwald. If $\alpha $ and $\tilde{F}$
are projectively related, then by Corollary 3.2 it results that $\tilde{F}$
is also generalized Berwald. So that, by (\ref{20}), (\ref{21}) and Lemma
4.1, the conditions (\ref{19})  reduce to $B^i=-P\eta ^i,$ for any $i=%
\overline{1,n},$ where $P=-\frac{\bar{\beta}}{2\tilde{F}|\beta |}\Gamma _{i%
\bar{j}}^kb_k\eta ^i\bar{\eta}^j.$

Conversely, if $\tilde{F}$ is generalized Berwald, then the first condition
from (\ref{19}) is identically satisfied and by (\ref{21}), $B^i=-\frac 1{%
\tilde{F}}\theta ^{*l}(\dot{\partial}_l\tilde{F})\eta ^i$ and $P=\frac 1{%
\tilde{F}}\theta ^{*i}(\dot{\partial}_i\tilde{F}).$ All these conditions
imply the projectiveness of the metrics $\alpha $ and $\tilde{F}.$ ii) is a
consequence of i), under assumptions of K\"{a}hler for the metrics $\alpha $
and $\tilde{F},$ respectively.
\end{proof}

\textbf{Example. }Let $\Delta =\left\{ (z,w)\in \mathbf{C}%
^2,\;|w|<|z|<1\right\} $ be the Hartogs triangle with the K\"{a}hler-purely
Hermitian metric
\begin{equation}
a_{i\overline{j}}=\frac{\partial ^2}{\partial z^i\partial \overline{z}^j}%
(\log \frac 1{\left( 1-|z|^2\right) \left( |z|^2-|w|^2\right) });\mbox{
}\alpha ^2(z,w;\eta ,\theta )=a_{i\overline{j}}\eta ^i\overline{\eta }^j,\;
\label{III.10}
\end{equation}
where $z,$ $w,\eta ,$ $\theta $ are the local coordinates $z^1,$ $z^2,$ $%
\eta ^1,$ $\eta ^2,$ respectively, and $|z^i|^2:=z^i\bar{z}^i,$ $z^i\in
\{z,w\},$ $\eta ^i\in \{\eta ,\theta \}.$ We choose
\begin{equation}
b_z=\frac w{|z|^2-|w|^2};\;b_w=-\frac z{|z|^2-|w|^2}.  \label{III.11}
\end{equation}
With these tools we have constructed in \cite{Al-Mu3} the complex Randers
metric $\tilde{F}=\alpha +|\beta |,$ where $\alpha (z,w,\eta ,\theta ):=%
\sqrt{a_{i\bar{j}}(z,w)\eta ^i\bar{\eta}^j}$ and $\beta (z,\eta
)=b_i(z,w)\eta ^i$. It is a complex Berwald metric, and so, by Theorem 4.3
ii), $\alpha $ and $\tilde{F}$ are projectively related.

Our second goal is to find when a complex Randers metric $\tilde{F}=\alpha
+|\beta |$ on a domain $D$ from $\mathbf{C}^{n}$ is projectively related to
the complex Euclidean metric $F$ on $D$.

For this, we make several assumptions. On the one hand we assume that $%
\tilde{F}$ is a complex Berwald metric. Thus, by Theorem 4.3, ii) we obtain
that $\alpha $ and $\tilde{F}$ are projectively related, $\alpha $ is
K\"{a}hler and $\tilde{G}^{i}=\stackrel{a}{G^{i}}.$ On the other hand, we
assume that $\alpha $ is projectively related to the Euclidean metric $F.$
Therefore, Theorem 3.7 implies that $\stackrel{a}{G^{i}}=\frac{1}{\alpha }%
\frac{\partial \alpha }{\partial z^{k}}\eta ^{k}\eta ^{i}$ . Under these
statements, we compute

$\frac 1{\tilde{F}}\frac{\partial \tilde{F}}{\partial z^k}\eta ^k\eta
^i=\frac 1{\tilde{F}}\frac{\partial \alpha }{\partial z^k}\eta ^k\eta
^i+\frac 1{\tilde{F}}\frac{\partial |\beta |}{\partial z^k}\eta ^k\eta ^i$

$=\frac 1{\tilde{F}}\frac{\partial \alpha }{\partial z^k}\eta ^k\eta
^i+\frac 1{2|\beta |\tilde{F}}\left( (\delta _k|\beta |)\eta ^k+2\bar{\beta}%
\stackrel{a}{G^l}b_l\right) \eta ^i$

$=\frac \alpha {\tilde{F}}\stackrel{a}{G^i}+\frac{|\beta |}{\tilde{F}}\frac
1\alpha \frac{\partial \alpha }{\partial z^k}\eta ^k\eta ^i=\stackrel{a}{G^i}%
.$ Thus, $\tilde{G}^i=\frac 1{\tilde{F}}\frac{\partial \tilde{F}}{\partial
z^k}\eta ^k\eta ^i,$ for any $i=\overline{1,n},$ which together with the
Berwald assumption for $\tilde{F},$ give that $\tilde{F}$ is projectively
related to the complex Euclidean metric $F$.

Conversely, by Theorem 3.7 it results that $F$ and $\tilde{F}$ are
projectively related if and only if the complex Randers metric $\tilde{F}$
is weakly K\"{a}hler and $\tilde{G}^i=\frac 1{\tilde{F}}\frac{\partial
\tilde{F}}{\partial z^k}\eta ^k\eta ^i$, for any $i=\overline{1,n}$. These
induce the generalized Berwald property for $\tilde{F}$ and by Theorem 4.2, $%
\tilde{F}$ is a complex Berwald metric. Now, taking into account Theorem
4.3, ii) it results that $\tilde{F}$ and $\alpha $ are projectively related,
$\alpha $ is K\"{a}hler and $\tilde{G}^i=\stackrel{a}{G^i}.$

So, we obtain
\begin{equation}
\stackrel{a}{G^i}=\frac 1{\tilde{F}}\left( \frac{\partial \alpha }{\partial
z^k}\eta ^k+\frac{\bar{\beta}}{|\beta |}\stackrel{a}{G^l}b_l\right) \eta ^i.
\label{22}
\end{equation}
The contraction with $b_i$ of (\ref{22}) gives $\stackrel{a}{G^i}b_i=\frac
\beta \alpha \frac{\partial \alpha }{\partial z^k}\eta ^k,$ which
substituted into (\ref{22}) yields $\stackrel{a}{G^i}=\frac 1\alpha \frac{%
\partial \alpha }{\partial z^k}\eta ^k\eta ^i$, i.e. $\alpha $ is
projectively related to the Euclidean metric $F.$

Therefore, the following theorem is proved

\begin{theorem}
Let $\tilde{F}=\alpha +|\beta |$ be a complex Randers metric on a domain $D$
from $\mathbf{C}^{n}$ and $F$ the complex Euclidean metric on $D.$ Then, $F$
and $\tilde{F}$ are projectively related if and only if $\alpha $ is
projectively related to the Euclidean metric $F$ and, $\tilde{F}$ is a
complex Berwald metric.
\end{theorem}

\textbf{Acknowledgment:} The first author is supported by the Sectorial
Operational Program Human Resources Development (SOP HRD), financed from the
European Social Fund and by Romanian Government under the Project number
POSDRU/89/1.5/S/59323.

\begin{flushleft}

Transilvania Univ.,
Faculty of Mathematics and Informatics

Iuliu Maniu 50, Bra\c{s}ov 500091, ROMANIA

e-mail: nicoleta.aldea@lycos.com

e-mail: gh.munteanu@unitbv.ro
\end{flushleft}

\end{document}